\documentclass[graybox,nospthms,leqno]{svmult}

\usepackage{mathptmx}       
\usepackage{helvet}         
\usepackage{courier}        
\usepackage{type1cm}        
\usepackage{makeidx}         
\usepackage{graphicx}        
\usepackage{multicol}        
\usepackage[bottom]{footmisc}

\usepackage{amsmath,amsthm}

\newcommand{\M}{\mathrm{M}}
\newcommand{\Orth}{\mathrm{O}}
\newcommand{\Sp}{\mathrm{Sp}}
\newcommand{\OM}{\mathrm{OM}}
\newcommand{\GO}{\mathrm{GO}}
\newcommand{\GSp}{\mathrm{GSp}}
\newcommand{\SpM}{\mathrm{SpM}}

\newcommand{\GL}{\mathrm{GL}}
\newcommand{\SL}{\mathrm{SL}}

\newcommand{\Hom}{\mathrm{Hom}}
\newcommand{\End}{\mathrm{End}}
\newcommand{\ind}{\mathrm{ind}}

\newcommand{\bil}[2]{\langle #1, #2 \rangle}
\newcommand{\transpose}{\mathsf{T}}
\newcommand{\Oh}{\mathcal{O}}

\swapnumbers 
\newtheorem{thm}[subsection]{Theorem}
\newtheorem{lem}[subsection]{Lemma}
\newtheorem{prop}[subsection]{Proposition}
\newtheorem{cor}[subsection]{Corollary}
\theoremstyle{definition}
\newtheorem{rmk}[subsection]{Remark}

\newtheorem{subsec}[subsection]{}

\begin{document}

\title*{Representations of reductive normal algebraic monoids}
\author{Stephen Doty}
\institute{
Loyola University Chicago, Department of
  Mathematics and Statistics, Chicago IL 60660 USA,
  \email{doty@math.luc.edu} 
}
%
%

\maketitle

\centerline{\emph{Dedicated to Lex Renner and Mohan Putcha}}

\bigskip

\abstract{The rational representation theory of a reductive normal
  algebraic monoid (with one-dimensional center) forms a highest
  weight category, in the sense of Cline, Parshall, and Scott. This is
  a fundamental fact about the representation theory of reductive
  normal algebraic monoids. We survey how this result was obtained,
  and treat some natural examples coming from classical groups.}

\section*{Introduction}\label{sec:intro}\noindent
Let $M$ be an affine algebraic monoid over an algebraically closed
field $K$.  See \cite{Putcha, Solomon, Renner:book} for general
surveys and background on algebraic monoids.  Assuming that $M$ is
reductive (its group $G$ of units is a reductive group) what can be
said about the representation theory of $M$ over $K$?

Recall that any affine algebraic group is smooth and hence normal (as
a variety). The normality of the algebraic group plays a significant
role in its representation theory, for instance in the proof of
Chevalley's theorem classifying the irreducible representations. Thus
it seems reasonable in trying to extend (rational) representation
theory from reductive groups to reductive monoids to look first at the
case when the monoid $M$ is normal. Furthermore, even in cases where a
given reductive algebraic monoid is not normal, one may always pass to
its normalization, which should be closely related to the original
object.

Renner \cite{Renner:classification} has obtained a classification
theorem for reductive normal algebraic monoids under the additional
assumptions that the center $Z(M)$ is 1-dimensional and that $M$ has a
zero element. Renner's classification theorem depends on an algebraic
monoid version of Chevalley's big cell, which holds for any reductive
affine algebraic monoid (with no assumptions about its center or a
zero). As a corollary of its construction, Renner derives a very
useful “extension principle” \cite[(4.5)]{Renner:classification} which
is a key ingredient in the analysis.

\section{Reductive normal algebraic monoids}\noindent
Let $M$ be a \emph{linear algebraic monoid} over an algebraically
closed field $K$. In other words, $M$ is a monoid (with unit element
$1 \in M$) which is also an affine algebraic variety over $K$, such
that the multiplication map $\mu: M \times M \to M$ is a morphism of
varieties. We assume that $M$ is irreducible as a variety. Hence the
unit group $G = M^\times$ (the subgroup of invertible elements of $M$)
is a connected linear algebraic group over $K$ and $G$ is Zariski
dense in $M$.

\begin{subsec}
Associated with $M$ is its affine coordinate algebra $K[M]$, the ring
of regular functions on $M$. There exist $K$-algebra homomorphisms
\[
\Delta: K[M] \to K[M] \otimes_K K[M], \quad \varepsilon: K[M] \to K
\]
called comultiplication and counit, respectively. For a given $f \in
K[M]$, we have $\varepsilon(f) = f(1)$; furthermore, if $\Delta(f) =
\sum_{i=1}^r f_i \otimes f'_i$ then $f(m_1m_2) = \sum_{i=1}^r
f_i(m_1)f'_i(m_2)$, for all $m_1, m_2 \in M$. The maps $\Delta,
\varepsilon$ make $K[M]$ into a bialgebra over $K$. This means that
they satisfy the bialgebra axioms:

(1) \quad $(id \otimes \Delta) \circ \Delta = (\Delta \otimes id)
\circ \Delta$,

(2) \quad $(\varepsilon\ \overline{\otimes}\ id) \circ \Delta = id =
(id\ \overline{\otimes}\ \varepsilon) \circ \Delta$

\noindent
where $\varphi\ \overline{\otimes}\ \varphi'$ denotes the map $a
\otimes a' \mapsto \varphi(a) \varphi'(a')$.

We note that the commutative bialgebra $(K[M], \Delta, \varepsilon)$
determines $M$, as the set $\Hom_{K\mathrm{-alg}}(K[M],K)$ of
$K$-algebra homomorphisms from $K[M]$ into $K$. The multiplication on
this set is defined by $\varphi \cdot \varphi' = (\varphi\
\overline{\otimes}\ \varphi') \circ \Delta$ and the identity element
is just the counit $\varepsilon$. One easily verifies that this
reconstructs $M$ from its coordinate bialgebra $K[M]$.

More generally, given any commutative bialgebra $(A,\Delta,
\varepsilon)$ over $K$, one defines on the set $M(A) =
\Hom_{K\mathrm{-alg}}(A,K)$ an algebraic monoid structure with
multiplication $\mu(\varphi, \varphi') = \varphi \cdot \varphi' =
(\varphi\ \overline{\otimes}\ \varphi') \circ \Delta$. This gives a
functor
\[
\{\text{commutative bialgebras over } K\} \to \{\text{algebraic
  monoids over } K\}
\]
which is quasi-inverse to the functor $M
\mapsto K[M]$. Thus the two categories are antiequivalent.
\end{subsec}

\begin{subsec}
Since $G$ is dense in $M$, the restriction map $K[M] \to K[G]$ (given
by $f \mapsto f_{|G}$) is injective, so we may identify $K[M]$ with a
subbialgebra of the Hopf algebra $K[G]$ of regular functions on $G$.
\end{subsec}

\begin{subsec}
Assume that $M$ is reductive; i.e., its unit group $G = M^\times$ is
reductive as an algebraic group. Fix a maximal torus $T$ in $G$. (Up
to conjugation $T$ is unique.) Let $X(T) = \Hom(T,K^\times)$ be the
character group of $T$; this is the abelian group of morphisms from
$T$ into the multiplicative group $K^\times$ of $K$. Let $X^\vee(T) =
\Hom(K^\times, T)$ be the abelian group of cocharacters into $T$. Let
$R \subset X(T)$ be the root system for the pair $(G,T)$ and $R^\vee
\subset X^\vee(T)$ the system of coroots.  According to the
classification of reductive algebraic groups, the reductive group $G$
is uniquely determined up to isomorphism by its root datum $(X(T), R,
X^\vee(T), R^\vee)$.
\end{subsec}

\begin{subsec}
We now add the assumption that $M$ is \emph{normal} as a variety. 
Let $D = \overline{T}$ be the Zariski closure of $T$ in $M$. Then $T
\subset D$ is an affine torus embedding. Let $X(D) = \Hom(D,K)$ be the
monoid of algebraic monoid homomorphisms from $D$ into $K$. The
restriction $\chi_{|T}$ of any $\chi \in X(D)$ is an element of
$X(T)$, so restriction defines a homomorphism $X(D) \to X(T)$. Since
$T$ is dense in $D$, this map is injective, and thus we may identify
$X(D)$ with a submonoid of $X(T)$.  Renner has shown that the
additional datum $X(D)$ is all that is needed to determine $M$ up to
isomorphism, under the additional hypotheses (probably unnecessary)
that the center 
\[
Z(M) = \{z \in M: zm=mz, \text{ for all } m \in M\}
\]
is 1-dimensional and that $M$ has a zero element. (One can always add
a zero formally, so the last requirement is insubstantial.)

It turns out that the set $X(D)$ also determines the rational
representation theory of the reductive normal algebraic monoid $M$, in
a sense made precise in Section 3.
\end{subsec}

\begin{subsec}
Note that it is easy to construct reductive algebraic monoids. Start
with a rational representation $\rho: G \to \End_K(V)$ of a reductive
group $G$ in some vector space $V$ with $\dim_K V = n < \infty$. The
image $\rho(G)$ is a reductive affine algebraic subgroup of $\End_K(V)
\simeq \M_n(K)$, the monoid of all $n \times n$ matrices under
ordinary matrix multiplication. Desiring our monoid to have a center
of at least dimension 1, we include the scalars $K^\times$ as scalar
matrices, defining $G_0$ to be the subgroup of $\End_K(V)$ generated
by $\rho(G)$ and $K^\times$. Now we set $M = \overline{G_0}$, the
Zariski closure of $G_0$ in $\End_K(V) \simeq M_n(K)$. This is a
reductive algebraic monoid.

For example, if $G = \SL_n(K)$ and $V$ is its natural representation
then $G_0 \simeq \GL_n(K)$ and $M = \M_n(K)$. (In general, to obtain a
monoid $M$ closely related to the starting group $G$, one should pick
$V$ to be a faithful representation.) There is no guarantee that this
procedure will always produce a normal reductive monoid, but if not
then one can always pass to its normalization.
\end{subsec}

\section{Examples: symplectic and orthogonal monoids}\noindent
The paper \cite{Doty:PRAMSACT} considered some more substantial
examples of reductive algebraic monoids coming from other classical
groups.  Let $V = K^n$ with its standard basis $\{e_1, \dots,
e_n\}$. Put $i^\prime = n+1-i$ for any $i = 1, \dots, n$.

\begin{subsec}\textbf{The orthogonal monoid.} 
Assume the characteristic of $K$ is
not 2. Define a symmetric nondegenerate bilinear form $\bil{\ }{\ }$
on $V$ by putting
\begin{equation*}
  \bil{e_i}{e_j} = \delta_{i, j'} \qquad \text{for any } 1 \le i,j \le
    n.  \tag{a}
\end{equation*}
Here $\delta$ is Kronecker's delta function.  Let $J$ be the matrix of
$\bil{\ }{\ }$ with respect to the basis $\{e_1, \dots, e_n\}$. Then
the orthogonal group $\Orth(V)$ is the group of linear operators $f
\in \End_K(V)$ preserving the form:
\begin{equation*}
\Orth(V) = \{ f \in \End_K(V): \bil{f(v)}{f(v')} = \bil{v}{v'},
\text{ all } v,v' \in V\}. \tag{b}
\end{equation*} 
Let $A$ be the matrix of $f$ with respect to the basis $\{e_1, \dots,
e_n\}$. Then we may identify $\Orth(V)$ with the matrix group
\begin{equation*}
\Orth_n(K) = \{ A \in \M_n(K): A^\transpose J A = J \}. \tag{c}
\end{equation*}
This is contained in the larger group $\GO_n(K)$, the group of
orthogonal similitudes (see e.g., \cite{involutions}) defined by
\begin{equation*}
\GO_n(K) = \{ A \in \M_n(K): A^\transpose J A = cJ, \text{ some } c
\in K^\times \}. \tag{d}
\end{equation*}
Note that $\GO_n(K)$ is generated by $\Orth_n(K)$ and $K^\times$. We
define the \emph{orthogonal monoid} $\OM_n(K)$ to be
\begin{equation*}
\OM_n(K) = \overline{\GO_n(K)}, \tag{e}
\end{equation*}
the Zariski closure in $\M_n(K)$. These monoids (for $n$ odd) were
studied by Grigor'ev \cite{Grigorev}. In \cite{Doty:PRAMSACT} the
following result was proved.
\end{subsec}

\begin{prop}\label{prop:2.2}
  The orthogonal monoid $\OM_n(K)$ is the set of all $A \in \M_n(K)$
  such that $A^\transpose J A = cJ = AJA^\transpose,$ for some $c
  \in K$.
\end{prop}

Note that the scalar $c\in K$ in the above is allowed to be zero, and
the ``extra'' condition $cJ = AJA^\transpose$ is necessary. If $c \ne
0$ then it is easy to see that $A^\transpose J A = cJ$ is equivalent
to $cJ = AJA^\transpose$, but when $c=0$ this equivalence fails. The
equivalence means that we could just as well have defined $\GO_n(K)$ by
\[
\GO_n(K) = \{ A \in \M_n(K): A^\transpose J A = cJ = AJA^\transpose,
\text{ some } c \in K^\times \}
\]
which is perhaps more suggestive for the description of $\OM_n(K)$
given above.

\begin{subsec}\textbf{The symplectic monoid.}
Assume that $n = \dim_K V$ is even, say $n = 2m$. Define an
antisymmetric nondegenerate bilinear form $\bil{\ }{\ }$ on $V$ by
putting
\begin{equation*}
\bil{e_i}{e_j} = \varepsilon_i\delta_{i, j'} \qquad \text{for any } 1
  \le i,j \le n. \tag{a}
\end{equation*}
where $\varepsilon_i$ is 1 if $i\le m$ and $-1$ otherwise.  Let $J$ be
the matrix of $\bil{\ }{\ }$ with respect to the basis $\{e_1, \dots,
e_n\}$. Then the symplectic group $\Sp(V)$ is the group of linear
operators $f \in \End_K(V)$ preserving the bilinear form:
\begin{equation*}
\Sp(V) = \{ f \in \End_K(V): \bil{f(v)}{f(v')} = \bil{v}{v'},
\text{ all } v,v' \in V\}. \tag{b}
\end{equation*} 
Let $A$ be the matrix of $f$ with respect to the basis $\{e_1, \dots,
e_n\}$. Then we may identify $\Sp(V)$ with the matrix group
\begin{equation*}
\Sp_n(K) = \{ A \in \M_n(K): A^\transpose J A = J \}. \tag{c}
\end{equation*}
This is contained in the larger group $\GSp_n(K)$, the group of
symplectic similitudes, defined by 
\begin{equation*}
\GSp_n(K) = \{ A \in \M_n(K): A^\transpose J A = cJ, \text{ some } c
\in K^\times \}. \tag{d}
\end{equation*}
Note that $\GSp_n(K)$ is generated by $\Sp_n(K)$ and $K^\times$. As in
the orthogonal case, we could just as well have defined $\GSp_n(K)$ by
\[
\GSp_n(K) = \{ A \in \M_n(K): A^\transpose J A = cJ = AJA^\transpose,
\text{ some } c \in K^\times \} 
\]
thanks to the equivalence of the conditions $A^\transpose J A = cJ$
and $cJ = AJA^\transpose$ in case $c\ne 0$.  We define the
\emph{symplectic monoid} $\SpM_n(K)$ to be
\begin{equation*}
\SpM_n(K) = \overline{\GSp_n(K)}, \tag{e}
\end{equation*}
the Zariski closure in $\M_n(K)$.  In \cite{Doty:PRAMSACT} the
following was proved.
\end{subsec}

\begin{prop}\label{prop:2.4}
  The symplectic monoid $\SpM_n(K)$ is the set of all $A \in \M_n(K)$
  such that $A^\transpose J A = cJ = AJA^\transpose,$ for some $c \in
  K$.
\end{prop}

Note that the scalar $c\in K$ in the above is allowed to be zero, and
the condition $cJ = AJA^\transpose$ is necessary, just as it was in the
orthogonal case.

\begin{subsec}\textbf{Sketch of proof.} 
I want to briefly sketch the ideas involved in the proof of
Propositions \ref{prop:2.2} and \ref{prop:2.4}.  Full details are
available in \cite{Doty:PRAMSACT}. The method of proof works for any
infinite field $K$ (except that characteristic 2 is excluded in the
orthogonal case). We continue to assume that $n=2m$ is even in the
symplectic case.

Let $G = \GO_n(K)$ or $\GSp_n(K)$ and let $M = \OM_n(K)$ or
$\SpM_n(K)$, respectively. Let $T$ be the maximal torus of diagonal
elements of $G$. Then we have inclusions
\begin{equation*}
\overline{T} \subset \overline{G} \subset M \tag{a}
\end{equation*}
and we desire to prove that the latter inclusion is actually an
equality. To accomplish this, we consider the action of $G \times G$
on $M$ given by $(g,h) \cdot m = gm h^{-1}$. Suppose that we can show
that every $G \times G$ orbit is of the form $GaG$, for some $a \in
\overline{T}$. Then it follows that
\begin{equation*}
  \textstyle  M = \bigcup_{a \in \overline{T}} GaG \subset
  \overline{G} \tag{b}
\end{equation*}
and this gives the opposite inclusion that proves Propositions
\ref{prop:2.2} and \ref{prop:2.4}. In fact, as it turned out, the
distinct $a \in \overline{T}$ in the above decomposition can be taken
to be certain idempotents in $\overline{T}$.

This suggests the program that was carried out in
\cite{Doty:PRAMSACT}, which in the end leads to additional structural
information on $M$:

(i) Classify all idempotents in $\overline{T}$.

(ii) Obtain an explicit description of $\overline{T}$. 

(iii) Determine the $G \times G$ orbits in $M$. 

\noindent
Part (i) is easy. For part (ii) one exploits the action of $T$ on
$\overline{T}$ by left multiplication and determines the orbits of
that action. Part (iii) involves developing orthogonal and symplectic
versions of classical Gaussian elimination.
\end{subsec}

\begin{subsec}\textbf{The normality question.} 
It is clear from the equalities
in Propositions \ref{prop:2.2} and \ref{prop:2.4} that $\OM_n(K)$ and
$\SpM_n(K)$ both have one-dimensional centers and contain zero. What
is not clear, and not addressed in \cite{Doty:PRAMSACT}, is whether or
not they are normal as algebraic varieties.

This question was recently settled in \cite{Doty-Hu}, where it is
shown that $\SpM_n(K)$ is always normal, while $\OM_n(K)$ is normal
only in case $n$ is even. More precisely, it is shown in
\cite{Doty-Hu} that when $n=2m$ is even, $\OM_n^+(K)$ and $\OM_n^-(K)$
are both normal varieties. Here
\begin{equation*}
\OM_n(K) = \OM_n^+(K) \cup \OM_n^-(K) \tag{a}
\end{equation*}
is the decomposition into irreducible components, where $\OM_n^+(K)$
is the component containing the unit element $1$.
\end{subsec}

\section{Representation theory}\noindent
From now on we assume that $M$ is an arbitrary reductive normal
algebraic monoid, with unit group $G = M^\times$. We wish to describe
some results of \cite{Doty:RNAM}. The main result is that the category
of rational $M$-modules is a highest weight category in the sense of
Cline--Parshall--Scott \cite{CPS}.

\begin{subsec}
We work with a fixed maximal torus $T \subset G$, and set $D =
\overline{T}$.  We assume that $\dim Z(M) = 1$ and $0 \in M$.  Recall
that restriction induces an injection $X(D) \to X(T)$, so we may
identify $X(D)$ with a submonoid of $X(T)$. We fix a Borel subgroup
$B$ with $T \subset B \subset G$ and let the set $R^-$ of
\emph{negative} roots be defined by the pair $(B,T)$. We set $R^+ =
-(R^-)$, the set of positive roots. We have $R = R^+ \cup R^-$.  Let
\[
X(T)^+ = \{\lambda \in X(T): (\alpha^\vee, \lambda) \ge 0, \text{ for
  all } \alpha\in R^+\}
\]
be the usual set of dominant weights. We define
\[
X(D)^+ = X(T)^+ \cap X(D).
\]
\end{subsec}

\begin{subsec}
By a (left) rational $M$-module we mean a linear action $M \times V
\to V$ such that the coefficient functions $M \to K$ of the action are
all in $K[M]$. This is the same as having a (right) $K[M]$-comodule
structure on $V$. This means that we have a comodule structure map
\begin{equation*}
\Delta_V : V \to V \otimes_K K[M]. \tag{a}
\end{equation*}
Since $K[M] \subset K[G]$ the map $\Delta_V$ induces a corresponding
map $V \to V \otimes_K K[G]$ making $V$ into a $K[G]$-comodule; i.e.,
a rational $G$-module. Thus, rational $M$-modules may also be regarded
as rational $G$-modules.  Any rational $M$-module is semisimple when
regarded as a rational $D$-module, with corresponding weight space
decomposition
\begin{equation*}
  \textstyle V = \bigoplus_{\lambda \in X(D)} V_\lambda \tag{b}
\end{equation*}
where $V_\lambda = \{ v\in V: d \cdot v = \lambda(d)\, v, \text{ all }
d \in D\}$. 

Recall that any rational $G$-module $V$ is semisimple when regarded as
a rational $T$-module, with corresponding weight space decomposition
\begin{equation*}
  \textstyle V = \bigoplus_{\lambda \in X(T)} V_\lambda \tag{c}
\end{equation*}
where $V_\lambda = \{ v\in V: t \cdot v = \lambda(t)\, v, \text{ all }
t \in T\}$. If $V$ is a rational $M$-module then the weight spaces
relative to $T$ are the same as the weight spaces relative to $D$. 
So the weights of a rational $M$-module all belong to
$X(D)$. Conversely, we have the following.
\end{subsec}

\begin{lem}
  If $V$ is a rational $G$-module such that $$\{\lambda \in
  X(T): V_\lambda \ne 0\} \subset X(D)$$ then $V$ extends uniquely to
  a rational $M$-module.
\end{lem}

This is proved as an application of Renner's extension principle,
which is a version of Chevalley's big cell construction for algebraic
monoids.

\begin{rmk}
  A special case of the lemma (for the case $M = \M_n(K)$) can be
  found in \cite{FS}.
\end{rmk}

\begin{subsec}
Next one needs a notion of induction for algebraic monoids, i.e., a
left adjoint to restriction. The usual definition of induced module
for algebraic groups does not work for algebraic monoids. Instead, we
use the following definition. Let $V$ be a rational $L$-module where
$L$ is an algebraic submonoid of $M$. We define $\ind_L^M V$ by
\[
\ind_L^M V = \{ f \in \Hom(M,V): f(lm) = l\cdot f(m), \text{ all } l
\in L, m \in M\}.
\]
This is viewed as a rational $M$-module via right translation. One can
check that in case $L,M$ are algebraic groups then this is isomorphic
to the usual induced module. 

It is well known that the Borel subgroup $B$ has the decomposition $B
= TU$, where $U$ is its unipotent radical. Given a character $\lambda
\in X(T)$ one regards $K$ as a rational $T$-module via $\lambda$; this
is often denoted by $K_\lambda$. One extends $K_\lambda$ to a rational
$B$-module by letting $U$ act trivially. Similarly, we have the
decomposition $\overline{B} = DU$. If $\lambda \in X(D)$ then we have
$K_\lambda$ as above, and again we may regard this as a rational
$\overline{B}$-module by letting $U$ act trivially. 
\end{subsec}

Now we can formulate the classification of simple rational
$M$-modules.

\begin{thm}
  Let $M$ be a reductive normal algebraic monoid. Let $\lambda \in
  X(D)$ and let $K_\lambda$ be the rational $\overline{B}$-module as
  above. Then

  {\rm(a)} $\ind_{\overline{B}}^M K_\lambda \ne 0$ if and only if
  $\lambda \in X(D)^+$.  

  {\rm(b)} If $\ind_{\overline{B}}^M K_\lambda \ne 0$ then its socle
  is a simple rational $M$-module (denoted by $L(\lambda)$).

  {\rm(c)} The set of $L(\lambda)$ with $\lambda \in X(D)^+$ gives a
  complete set of isomorphism classes of simple rational $M$-modules.
\end{thm}

Let $\lambda \in X(T)$. Let $\nabla(\lambda) = \ind_{B}^G
K_\lambda$. It is well known that $\nabla(\lambda) \ne 0$ if and only
if $\lambda \in X(T)^+$.  The following is a key fact.

\begin{lem}
  If $\lambda \in X(D)^+$ then $\ind_{\overline{B}}^M K_\lambda =
  \nabla(\lambda) = \ind_B^G K_\lambda$.
\end{lem}

\begin{subsec}
Now we consider truncation. Let $\pi \subset X(T)^+$.  Given a
rational $G$-module $V$, let $\Oh_\pi V$ be the unique largest
rational submodule of $V$ with the property that the highest weights
of all its composition factors belong to $\pi$. The (left exact)
truncation functor $\Oh_\pi$ was defined by Donkin \cite{Donkin:SA1}.

Recall that $X(T)$ is partially ordered by $\lambda \le \mu$ if $\mu -
\lambda$ can be written as a sum of positive roots; this is sometimes
called the dominance order.  A subset $\pi$ of $X(T)^+$ is said to be
\emph{saturated} if it is predecessor closed under the dominance order
on $X(T)$. In other words, $\pi$ is saturated if for any $\mu \in \pi$
and any $\lambda \in X(T)^+$, $\lambda \le \mu$ implies that $\lambda
\in \pi$.
\end{subsec}

In order to show that the category of rational $M$-modules is a
highest weight category, we are going to take $\pi = X(D)^+$.  We need
the following observation.

\begin{lem}
  The set $\pi = X(D)^+$ is a saturated subset of $X(T)^+$. 
\end{lem}

For $\lambda \in X(T)^+$, let $I(\lambda)$ be the injective envelope
of $L(\lambda)$ in the category of rational $G$-modules. For $\lambda
\in X(D)^+$ let $Q(\lambda)$ be the injective envelope of $L(\lambda)$
in the category of rational $M$-modules. The following records the
effect of truncation on various classes of rational $G$-modules.

\begin{thm}
  Let $\pi = X(D)^+$. For any $\lambda \in X(T)^+$ we have the
  following:

  {\rm(a)}\quad $\Oh_\pi
  \nabla(\lambda) =
  \begin{cases}
    \nabla(\lambda) & \text{ if } \lambda \in \pi \\
    0 & \text{ otherwise.}
  \end{cases}$

  {\rm(b)}\quad $\Oh_\pi I(\lambda) =
  \begin{cases}
    Q(\lambda) & \text{ if } \lambda \in \pi \\
    0 & \text{ otherwise.}
  \end{cases}$

  {\rm(c)}\quad $\Oh_\pi K[G] = K[M]$.
\end{thm}

Note that $K[M]$ is regarded as a rational $M$-module via right
translation.  A $\nabla$-filtration for a rational $G$-module $V$ is
an ascending series $$0 = V_0 \subset V_1 \subset \cdots \subset
V_{r-1} \subset V_r = V$$ of rational submodules such that for each $j
= 1, \dots, r$, the quotient $V_j/V_{j-1}$ is isomorphic to some
$\nabla(\lambda_j)$. Whenever $V$ is a rational $G$-module with a
$\nabla$-filtration, let $(V: \nabla(\lambda))$ be the number of
$\lambda_j$ for which $\lambda = \lambda_j$. This number is
independent of the filtration. 

The proof of the above theorem, which relies on results of
\cite{Donkin:filt}, also shows the following facts.

\begin{cor}
  {\rm(a)} Let $\lambda \in \pi = X(D)^+$. The module $Q(\lambda)$ has
  a $\nabla$-filtration. Furthermore, it satisfies the reciprocity
  property
  \[
  (Q(\lambda): \nabla(\mu)) = [\nabla(\mu): L(\lambda)]
  \]
  for any $\mu \in X(D)^+$, where $[V:L]$ stands for the multiplicity
  of a simple module $L$ in a composition series of $V$.

  {\rm(b)} The module $K[M]$ has a $\nabla$-filtration. Moreover,
  $(K[M]: \nabla(\lambda)) = \dim_K \nabla(\lambda)$ for each $\lambda
  \in X(D)^+$.  
\end{cor}

\begin{subsec}
From these results one obtains the important fact that the category of
rational $M$-modules is a highest weight category, in the sense of
\cite{CPS}. In particular, one also sees that $\dim_K Q(\lambda)$ is
finite, for any $\lambda \in X(D)^+$. (In contrast, it is well known
that $\dim_K I(\lambda)$ is infinite.) 
\end{subsec}

\begin{subsec}
It is also shown in \cite{Doty:RNAM}, exploiting the assumption that
$Z(M)$ is one-dimensional, that the category of rational $M$-modules
splits into a direct sum of `homogeneous' subcategories each of which
is controlled by a finite saturated subset of $X(D)^+$. From the
results of \cite{CPS} it then follows that there is a finite
dimensional quasihereditary algebra in each homogeneous degree, whose
module category is precisely the homogeneous subcategory in that
degree. Details are given in \cite{Doty:RNAM}, where it is also shown
that the quasihereditary algebras in question are in fact generalized
Schur algebras in the sense of Donkin \cite{Donkin:SA1}.
\end{subsec}

\bibliographystyle{amsalpha}

\end{document}